\documentclass[a4paper,12pt]{article}
\usepackage{amscd}
\usepackage{amsmath,amsfonts,amssymb,amscd}
\usepackage{indentfirst,graphicx,epsfig}
\usepackage{graphicx,psfrag}
\input{epsf}

\newtheorem{thm}{Theorem}

\newtheorem{lem}{Lemma}

\def\qed{\hfill \nopagebreak\rule{5pt}{8pt}}
\def\pf{\noindent {\it Proof.} }

\title{\bf On the rainbow vertex-connection\footnote{Supported by NSFC and ``the
Fundamental Research Funds for the Central Universities". } }
\author{
\small  Xueliang Li, Yongtang Shi\\
\small Center for Combinatorics and LPMC-TJKLC \\
\small Nankai University, Tianjin 300071, China \\
\small Email: lxl@nankai.edu.cn,  shi@nankai.edu.cn
\date{}}
\begin{document}
\maketitle
\begin{abstract}
A vertex-colored graph is {\it rainbow vertex-connected} if any two
vertices are connected by a path whose internal vertices have
distinct colors, which was introduced by Krivelevich and Yuster. The
{\it rainbow vertex-connection} of a connected graph $G$, denoted by
$rvc(G)$, is the smallest number of colors that are needed in order
to make $G$ rainbow vertex-connected. Krivelevich and Yuster proved
that if $G$ is a graph of order $n$ with minimum degree $\delta$,
then $rvc(G)<11n/\delta$. In this paper, we show that $rvc(G)\leq
3n/(\delta+1)+5$ for $\delta\geq \sqrt{n-1}-1$ and $n\geq 290$,
while $rvc(G)\leq 4n/(\delta+1)+5$ for $16\leq \delta\leq
\sqrt{n-1}-2$ and $rvc(G)\leq 4n/(\delta+1)+C(\delta)$ for
$6\leq\delta\leq 15$, where
$C(\delta)=e^{\frac{3\log(\delta^3+2\delta^2+3)-3(\log 3-1)}
{\delta-3}}-2$. We also prove that $rvc(G)\leq 3n/4-2$ for
$\delta=3$, $rvc(G)\leq 3n/5-8/5$ for $\delta=4$ and $rvc(G)\leq
n/2-2$ for $\delta=5$. Moreover, an example shows that when
$\delta\geq \sqrt{n-1}-1$ and $\delta=3,4,5$, our bounds are
seen to be tight up to additive factors.\\
[2mm] {\bf Keywords:} rainbow vertex-connection; minimum degree;
$2$-step dominating set\\[2mm]
{\bf AMS subject classification 2010:} 05C15, 05C40
\end{abstract}

\section{Introduction}
All graphs considered in this paper are simple, finite and
undirected. We follow the notation and terminology of Bondy and
Murty \cite{BM}. An edge-colored graph is {\it rainbow connected} if
any two vertices are connected by a path whose edges have distinct
colors. Obviously, if $G$ is rainbow connected, then it is also
connected. This concept of rainbow connection in graphs was
introduced by Chartrand et al. in \cite{CJMZ}. The rainbow
connection number of a connected graph $G$, denoted by $rc(G)$, is
the smallest number of colors that are needed in order to make $G$
rainbow connected. Observe that $diam(G)\leq rc(G)\leq n-1$. It is
easy to verify that $rc(G)=1$ if and only if $G$ is a complete
graph, that $rc(G)=n-1$ if and only if $G$ is a tree. It was shown
that computing the rainbow connection number of an arbitrary graph
is NP-hard \cite{CFMY}. There are some approaches to study the
bounds of $rc(G)$, we refer to \cite{CLRTY,CDRV,KY,S}.

In \cite{KY}, Krivelevich and Yuster proposed the concept of rainbow
vertex-connection.  A vertex-colored graph is {\it rainbow
vertex-connected} if any two vertices are connected by a path whose
internal vertices have distinct colors. The {\it rainbow
vertex-connection} of a connected graph $G$, denoted by $rvc(G)$, is
the smallest number of colors that are needed in order to make $G$
rainbow vertex-connected. An easy observation is that if $G$ is of
order $n$ then $rvc(G)\leq n-2$ and $rvc(G)=0$ if and only if $G$ is
a complete graph. Notice that $rvc(G)\geq diam(G)-1$ with equality
if the diameter is $1$ or $2$. For rainbow connection and rainbow
vertex-connection, some examples are given to show that there is no
upper bound for one of parameters in terms of the other in
\cite{KY}.

Krivelevich and Yuster \cite{KY} proved that if $G$ is a graph with
$n$ vertices and minimum degree $\delta$, then $rvc(G)<11n/\delta$.
In this paper, we will improve this bound for given order $n$ and
minimum degree $\delta$. We will show that $rvc(G)\leq
3n/(\delta+1)+5$ for $\delta\geq \sqrt{n-1}-1$ and $n\geq 290$,
while $rvc(G)\leq 4n/(\delta+1)+5$ for $16\leq \delta\leq
\sqrt{n-1}-2$ and $rvc(G)\leq 4n/(\delta+1)+C(\delta)$ for
$6\leq\delta\leq 15$, where
$C(\delta)=e^{\frac{3\log(\delta^3+2\delta^2+3)-3(\log 3-1)}
{\delta-3}}-2$. We also prove that $rvc(G)\leq 3n/4-2$ for
$\delta=3$, $rvc(G)\leq 3n/5-8/5$ for $\delta=4$ and $rvc(G)\leq
n/2-2$ for $\delta=5$.

Moreover, an example shows that when $\delta\geq \sqrt{n-1}-1$ and
$\delta=3,4,5$, our bounds are seen to be tight up to additive
factors. To see this, we recall the graph constructed by Caro et al.
\cite{CJMZ}, which was used to interpret the upper bound of rainbow
connection $rc(G)$. A connected $n$-vertex graph $H$ are constructed
as follows. Take $m$ copies of complete graph $K_{\delta+1}$,
denoted $X_1, \ldots, X_m$ and label the vertices of $X_i$ with
$x_{i,1},\ldots,x_{i,\delta+1}$. Take two copies of $K_{\delta+2}$,
denoted $X_0$, $X_{m+1}$ and similarly label their vertices. Now,
connect $x_{i,2}$ with $x_{i+1,1}$ for $i = 0,\ldots, m$ with an
edge, and delete the edges $x_{i,1}x_{i,2}$ for i = $0,\ldots, m +
1$. Observe that the obtained graph $H$ has $n = (m + 2)(\delta + 1)
+ 2$ vertices, minimum degree $\delta$ and diameter $\frac {3n}
{\delta+1}-\frac{\delta+7}{\delta+1}$. Therefore, the upper bound of
$rvc(G)$ cannot be improved below $\frac {3n}
{\delta+1}-\frac{\delta+7}{\delta+1}-1=\frac {3n-6}{\delta+1}-2$.

\section{$rvc(G)$ and minimum degree}

Let $\mathcal{G}(n,\delta)$ be the class of simple connected
$n$-vertex graphs with minimum degree $\delta$. Let $\ell(n,\delta)$
be the maximum value of $m$ such that every $G\in
\mathcal{G}(n,\delta)$ has a spanning tree with at least $m$ leaves.
We can obtain an trivial upper bound for $rvc(G)$.

\begin{lem}
A connected graph $G$ of order $n$ with maximum degree $\Delta(G)$
has $rvc(G)\leq n-\ell(n,\delta)$ and $rvc(G)\leq n-\Delta(G)$.
\end{lem}
\pf It is obvious that $rvc(G)\leq n-\ell(n,\delta)$. For a
connected graph $G$ of order $n$ with maximum degree $\Delta(G)=k$,
a spanning tree $T$ of $G$ grown from a vertex $v$ with degree $k$
has maximum degree $\Delta(T)=k$. It deduces that $T$ has at least
$k$ leaves. Thus, $rvc(G)\leq n-\Delta(G)$ holds.\qed\\

Note that finding the maximum number of leaves in a spanning tree of
$G$ is NP-hard. Linial and Sturtevant (unpublished) proved that
$\ell(n,3)\geq n/4 + 2$. For $\delta = 4$, the optimal bound
$\ell(n,4)\geq 2/5n + 8/5$ is proved in \cite{GW} and in \cite{KW}.
In \cite{GW}, it is also proved that $\ell(n,5)\geq n/2 + 2$.
Indeed, Kleitman and West in \cite{KW} proved that
$\ell(n,\delta)\geq (1-b\ln \delta/\delta)n$ for large $\delta$,
where $b$ is any constant exceeding $2.5$. Hence, the following
theorem is obvious.

\begin{thm}
For a connected graph $G$ of order $n$ with minimum degree $\delta$,
$rvc(G)\leq 3n/4-2$ for $\delta=3$, $rvc(G)\leq 3n/5-8/5$ for
$\delta=4$ and $rvc(G)\leq n/2-2$ for $\delta=5$. For sufficiently
large $\delta$, $rvc(G)\leq (b \ln \delta)n/\delta$, where $b$ is
any constant exceeding $2.5$.\qed
\end{thm}

Motivated by the proof methods in \cite{KY}, we will prove the
following theorem by constructing a connected $\delta/3$-strong
$2$-step dominating set $S$ whose size is at most $3n/(\delta+1)-2$.

\begin{thm}\label{thm1}
A connected graph $G$ of order $n$ with minimum degree $\delta$ has
$rvc(G)\leq 3n/(\delta+1)+5$ for $\delta\geq \sqrt{n-1}-1$ and
$n\geq 290$, while $rvc(G)\leq 4n/(\delta+1)+5$ for $16\leq
\delta\leq \sqrt{n-1}-2$ and $rvc(G)\leq 4n/(\delta+1)+C(\delta)$
for $6\leq\delta\leq 15$, where
$C(\delta)=e^{\frac{3\log(\delta^3+2\delta^2+3)-3(\log 3-1)}
{\delta-3}}-2$.
\end{thm}

Now we state some lemmas that are needed to prove Theorem
\ref{thm1}. The first lemma is from \cite{KY}.

\begin{lem}\label{lem2}
If $G$ is a connected graph with minimum degree $\delta$, then it
has a connected spanning subgraph with minimum degree $\delta$ and
with less than $n(\delta+1/(\delta+1))$ edges.\qed
\end{lem}

A set of vertices $S$ of a graph $G$ is called a {\it $2$-step
dominating set}, if every vertex of $V(G)\setminus S$ has either a
neighbor in $S$ or a common neighbor with a vertex in $S$. In
\cite{KY}, the authors introduced a special type of $2$-step
dominating set. A $2$-step dominating set is {\it $k$-strong}, if
every vertex that is not dominated by it has at least $k$ neighbors
that are dominated by it. We call a $2$-step dominating set $S$ is
{\it connected}, if the subgraph induced by $S$ is connected.
Similarly, we can define the concept of {\it connected $k$-strong
$2$-step dominating set}.

\begin{lem}
If $G$ is a graph of order $n$ with minimum degree $\delta\geq 2$,
then $G$ has a connected $\delta/3$-strong $2$-step dominating set
$S$ whose size is at most $3n/(\delta+1)-2$.
\end{lem}
\pf For any $T\subseteq V(G)$, denote by $N^k(T)$ the set of all
vertices with distance exactly $k$ from $T$. We construct a
connected $\delta/3$-strong $2$-step dominating set $S$ as follows:

{\bf Procedure 1}. Initialize $S'=\{u\}$ for some $u\in V(G)$. As
long as $N^3(S')\neq \emptyset$, take a vertex $v\in N^3(S')$ and
add vertices $v,x_1,x_2$ to $S'$, where $vx_2x_1x_0$ is a shortest
path from $v$ to $S'$ and $x_0\in S'$.

{\bf Procedure 2}. Initialize $S=S'$ obtained from Procedure 1. As
long as there exists a vertex $v\in N^2(S)$ such that $|N(v)\cap
N^2(S)|\geq  2\delta/3 +1$, add vertices $v,y_1$ to $S$, where
$vy_1y_0$ is a shortest path from $v$ to $S$ and $y_0\in S$.

Clearly $S'$ remains connected after every iteration in Procedure 1.
Therefore, when Procedure 1 ends, $S'$ is a connected $2$-step
dominating set. Let $k_1$ be the number of iterations executed in
Procedure 1. Observe that when a new vertex from $N^3(S')$ is added
to $S'$, $|S'\cup N^1(S')|$ increases by at least $\delta +1$ in
each iteration. Thus, we have $k_1+1\leq \frac{|S'\cup
N^1(S')|}{\delta +1}=\frac{n-|N^2(S')|}{\delta +1}$. Furthermore,
$|S'|=3k_1+1\leq \frac{3\left(n-|N^2(S')|\right)}{\delta +1}-2$
since three more vertices are added in each iteration.

Notice that $S$ also remains connected after every iteration in
Procedure 2. When Procedure 2 ends, each $v\in N^2(S)$ has at most
$2\delta/3$ neighbors in $N^2(S)$, i.e., has at least $\delta/3$
neighbors in $N^1(S)$, so $S$ is a connected $\delta/3$-strong
$2$-step dominating set. Let $k_2$ be the number of iterations
executed in Procedure 2. Observe that when a new vertex from
$N^2(S)$ is added to $S$, $|N^2(S)|$ reduces by at least $2\delta/3
+2$ in each iteration. Thus, we have $k_2\leq
\frac{|N^2(S')|}{2\delta/3 +2}=\frac{3|N^2(S')|}{2\delta +6}$.
Furthermore,
$$|S|=|S'|+2k_2\leq \frac{3\left(n-|N^2(S')|\right)}{\delta +1}-2+\frac{6|N^2(S')|}{2\delta
+6}< \frac{3n}{\delta +1}-2.$$
\qed\\

Before proceeding, we first recall the Lov\'asz Local Lemma
\cite{AS}.\\

\noindent{\bf The Lov\'asz Local Lemma}~~ Let $A_1, A_2, \ldots,
A_n$ be the events in an arbitrary probability space. Suppose that
each event $A_i$ is mutually independent of a set of all the other
events $A_j$ but at most $d$, and that $P[A_i]\leq p$ for all $1\leq
i\leq n$. If $ep(d+1) < 1$, then $Pr[\bigwedge_{i=1}^n
\overline{A_i}]>0$.\qed\\

\noindent{\bf Proof of Theorem \ref{thm1}}~~ Suppose $G$ is a
connected graph of order $n$ with minimum degree $\delta$. By Lemma
\ref{lem2}, we may assume that $G$ has less than
$n(\delta+1/(\delta+1))$ edges. Let $S$ be a connected
$\delta/3$-strong $2$-step dominating set of $G$ with at most
$3n/(\delta+1)-2$ vertices.

Suppose $\delta\geq \sqrt{n-1}-1$. Observe that each vertex $v$ of
$N^1(S)$ has less than $(\delta+1)^2$ neighbors in $N^2(S)$, since
$(\delta+1)^2\geq n-1$ and $v$ has another neighbor in $S$. We
assign colors to graph $G$ as follows: distinct colors to each
vertex of $S$ and seven new colors to vertices of $N^1(S)$ such that
each vertex of $N^1(S)$ chooses its color randomly and independently
from all other vertices of $N^1(S)$. Hence, the total color we used
is at most
$$|S|+7\leq\frac{3n}{\delta +1}-2+7=\frac{3n}{\delta +1}+5.$$
For each vertex $u$ of $N^2(S)$, let $A_u$ be the event that all the
neighbors of $u$ in $N^1(S)$, denoted by $N_{1}(u)$, are assigned at
least two distinct colors. Now we will prove $Pr[A_u]>0$ for each
$u\in N^2(S)$. Notice that each vertex $u\in N^2(S)$ has at least
$\delta/3$ neighbors in $N^1(S)$ since $S$ is a connected
$\delta/3$-strong $2$-step dominating set of $G$. Therefore, we fix
a set $X(u)\subset N^1(S)$ of neighbors of $u$ with
$|X(u)|=\lceil\delta/3\rceil$. Let $B_u$ be the event that all of
the vertices in $X(u)$ receive the same color. Thus, $Pr[B_u]\leq
7^{-\lceil\delta/3\rceil+1}$. As each vertex of $N^1(S)$ has less
than $(\delta+1)^2$ neighbors in $N^2(S)$, we have that the event
$B_u$ is independent of all other events $B_v$ for $v\neq u$ but at
most $((\delta+1)^2-1)\lceil\delta/3\rceil$ of them. Since for
$\delta\geq \sqrt{n-1}-1$ and $n\geq 290$,
$$e\cdot 7^{~-\lceil\delta/3\rceil+1}
\left(((\delta+1)^2-1)\lceil\delta/3\rceil+1\right)<1,$$ by the
Lov\'asz Local Lemma, we have $Pr[A_u]>0$ for each $u\in N^2(S)$.
Therefore, for $N^1(S)$, there exists one coloring with seven colors
such that every vertex of $N^2(S)$ has at least two neighbors in
$N^1(S)$ colored differently. It remains to show that graph $G$ is
rainbow vertex-connected. Let $u,v$ be a pair of vertices such that
$u,v\in N^2(S)$. If $u$ and $v$ have a common neighbor in $N^1(S)$,
then we are done. Denote by $x_1,y_1$ and $x_2,y_2$, respectively,
the two neighbors of $u$ and $v$ in $N^1(S)$ such that the colors of
$x_i$ and $y_i$ are different for $i=1,2$. Without loss of
generality, suppose the colors of $x_1$ and $x_2$ are also
different. Indeed, there exists a required path between $u$ and $v$:
$ux_1w_1Pw_2x_2v$, where $w_i$ is the neighbor of $x_i$ in $S$ and
$P$ is the path connected $w_1$ and $w_2$ in $S$. All other cases of
$u,v$ can be checked easily.

From now on we assume $\delta\leq \sqrt{n-1}-2$. We partition
$N^1(S)$ to two subsets: $D_1=\{v\in N^1(S): \mbox{~$v$ has at least
$(\delta+1)^2$ neighbors in $N^2(S)$} \}$ and $D_2=N^1(S)\setminus
D_1$. Since $G$ has less than $n(\delta+1/(\delta+1))$ edges, we
have $|D_1|\leq n/(\delta+1)$. Denote by $L_1=\{v\in N^2(S):
\mbox{~$v$ has at least one neighbor in $D_1$} \}$ and
$L_2=N^2(S)\setminus L_1$.

Let $C(\delta)=5$ for $16\leq \delta\leq \sqrt{n-1}-2$ and
$C(\delta)=e^{\frac{3\log(\delta^3+2\delta^2+3)-3(\log 3-1)}
{\delta-3}}-2$ for $6\leq \delta\leq 15$. We assign colors to graph
$G$ as follows: distinct colors to each vertex of $S\cup D_1$ and
$C(\delta)+2$ new colors to vertices of $D_2$ such that each vertex
of $D_2$ chooses its color randomly and independently from all other
vertices of $D_2$. Hence, the total color we used is at most
$$|S|+|D_1|+C(\delta)+2\leq\frac{3n}{\delta +1}-2+\frac n {\delta+1}+
C(\delta)+2=\frac{4n}{\delta +1}+C(\delta).$$ For each vertex $u$ of
$L_2$, let $A_u$ be the event that all the neighbors of $u$ in $D_2$
are assigned at least two distinct colors. Now we will prove
$Pr[A_u]>0$ for each $u\in L_2$. Notice that each vertex $u\in L_2$
has at least $\delta/3$ neighbors in $D_1$. Therefore, we fix a set
$X(u)\subset D_1$ of neighbors of $u$ with
$|X(u)|=\lceil\delta/3\rceil$. Let $B_u$ be the event that all of
the vertices in $X(u)$ receive the same color. Thus, $Pr[B_u]\leq
\left( C(\delta)+2\right)^{-\lceil\delta/3\rceil+1}$. As each vertex
of $D_2$ has less than $(\delta+1)^2$ neighbors in $N^2(S)$, we have
that the event $B_u$ is independent of all other events $B_v$ for
$v\neq u$ but at most $((\delta+1)^2-1)\lceil\delta/3\rceil$ of
them. Since
$$e\cdot\left(C(\delta)+2\right)^{-\lceil\delta/3\rceil+1}
\left(((\delta+1)^2-1)\lceil\delta/3\rceil+1\right)<1,$$ by the
Lov\'asz Local Lemma, we have $Pr[A_u]>0$ for each $u\in L_2$.
Therefore, for $D_2$, there exists one coloring with $C(\delta)+2$
colors such that each vertex of $L_2$ has at least two neighbors in
$D_2$ colored differently.

Similarly, we can check that graph $G$ is also rainbow
vertex-connected in this case.

The proof is thus completed.

\qed\\


\begin{thebibliography}{1}

\bibitem{AS}  N. Alon and J.H. Spencer, The Probabilistic Method, 3rd ed,
Wiley, New York, 2008.

\bibitem{BM} J.A. Bondy and U.S.R. Murty, Graph Theory, GTM 244, Springer, 2008.

\bibitem{CLRTY}
Y. Caro, A. Lev, Y. Roditty, Z. Tuza and R. Yuster,  On rainbow
connection, {\it Electron J. Combin.} {\bf 15}(2008), R57.

\bibitem{CFMY}
S. Chakraborty, E. Fischer, A. Matsliah and R. Yuster, Hardness and
algorithms for rainbow connectivity, {\it J. Comb. Optim.}, in
press.

\bibitem{CDRV}
L. Chandran, A. Das, D. Rajendraprasad and N. Varma, Rainbow
connection number and connected dominating sets, arXiv:1010.2296v1.

\bibitem{CJMZ}
G. Chartrand, G.L. Johns, K.A. McKeon and P. Zhang, Rainbow
connection in graphs, {\it Math. Bohemica} {\bf 133}(2008), 85--98.

\bibitem{GW}
J.R. Griggs and M. Wu, Spanning trees in graphs with minimum degree
$4$ or $5$, {\it Discrete Math.} {\bf 104}(1992), 167--183.

\bibitem{KW}
D.J. Kleitman and D.B. West, Spanning trees with many leaves, {\it
SIAM J. Discrete Math.} {\bf 4}(1991), 99--106.

\bibitem{KY}
M. Krivelevich and R. Yuster, The rainbow connection of a graph is
(at most) reciprocal to its minimum degree, {\it J. Graph Theory}
{\bf 63}(2010), 185--191.

\bibitem{S}
I. Schiermeyer, Rainbow connection in graphs with minimum degree
three, IWOCA 2009, {\it LNCS} {\bf 5874}(2009), 432--437.


\end{thebibliography}
\end{document}